\documentclass[12pt]{article}
\usepackage{amsmath, amsthm, amssymb,amscd, verbatim}

\theoremstyle{plain}
\newtheorem{thm}{Theorem}[section]
\newtheorem{lem}[thm]{Lemma}
\newtheorem{cor}[thm]{Corollary}

\newtheorem{prop}[thm]{Proposition}

\theoremstyle{definition}
\newtheorem{defn}[thm]{Definition}

\theoremstyle{remark}
\newtheorem{rem}[thm]{Remark}

\newcommand{\nc}{\newcommand} 
\nc{\hb}{\mathbb} 
\nc{\M}{\mathcal} 
\nc{\mf}{\mathfrak}
\nc{\mbf}{\mathbf}
\nc{\DMO}{\DeclareMathOperator}

\newbox\noforkbox \newdimen\forklinewidth
\setbox0\hbox{$\textstyle\smile$}
\setbox1\hbox to \wd0{\hfil\vrule width \forklinewidth depth-2pt
 height 10pt \hfil}
\setbox\noforkbox\hbox{\lower 2pt\box1\lower 2pt\box0\relax}
\def\anchor{\mathop{\copy\noforkbox}\limits}
 
\setbox0\hbox{$\textstyle\smile$}
\setbox1\hbox to \wd0{\hfil{\sl /\/}\hfil}
\setbox2\hbox to \wd0{\hfil\vrule height 10pt depth -2pt width
               \forklinewidth\hfil}
\newbox\doesforkbox
\setbox\doesforkbox\hbox{\box1 \lower 2pt\box2\lower2pt\box0\relax}
\def\nanchor{\mathop{\copy\doesforkbox}\limits}
 
\nc{\cA}{{\M A}} \nc{\cB}{{\M B}} \nc{\cC}{{\M C}} \nc{\cD}{{\M D}}
\nc{\cE}{{\M E}} \nc{\cF}{{\M F}} \nc{\cG}{{\M G}} \nc{\cH}{{\M H}}
\nc{\cI}{{\M I}} \nc{\cJ}{{\M J}} \nc{\cK}{{\M K}} \nc{\cL}{{\M L}}
\nc{\cM}{{\M M}} \nc{\cN}{{\M N}} \nc{\cO}{{\M O}} \nc{\cP}{{\M P}}
\nc{\cQ}{{\M Q}} \nc{\cR}{{\M R}} \nc{\cS}{{\M S}} \nc{\cT}{{\M T}}
\nc{\cU}{{\M U}} \nc{\cV}{{\M V}} \nc{\cW}{{\M W}} \nc{\cX}{{\M X}}
\nc{\cY}{{\M Y}} \nc{\cZ}{{\M Z}}
\nc{\Aa}{{\hb A}} \nc{\Cc}{{\hb C}} \nc{\Gg}{{\hb G}}
\nc{\Nn}{{\hb N}} \nc{\Pp}{{\hb P}} 
\nc{\Qq}{{\hb Q}} \nc{\Rr}{{\hb R}} \nc{\Zz}{{\hb Z}}
\nc{\mfa}{{\mf a}} \nc{\mfb}{{\mf b}} \nc{\mfk}{{\mf k}}
\nc{\mfm}{{\mf m}} \nc{\mfp}{{\mf p}} \nc{\mfq}{{\mf q}}
\nc{\mfr}{{\mf r}}
\nc{\fP}{{\mf P}}
\DMO*{\trdeg}{td}
\DMO*{\spec}{Spec}
\DMO*{\fork}{\nanchor}
\DMO*{\dnf}{\anchor}
\DMO{\RU}{RU}
\DMO{\RM}{RM}
\DMO{\RC}{RC}
\DMO{\Hu}{Hull}
\DMO{\leg}{length}
\DMO{\area}{area}
\DMO{\dia}{diameter}
\DMO{\iso}{ISO}
\DMO{\dis}{dist}
\DMO{\grad}{grad}
\DMO{\vol}{volume}
\DMO{\gra}{grad}
\DMO{\hd}{nbhd}
\DMO{\dv}{div}
\DMO{\Psl}{PSL}
\DMO{\Ric}{Ricci}
\DMO{\Hom}{Hom}
\DMO{\J}{Jac}
\DMO{\Id}{I}
\DMO{\inj}{injrad} 
\DMO{\supp}{supp}
\DMO{\Lip}{Lip}
\DMO{\SO}{SO}
\DMO{\local}{loc}
\nc{\Mb}{\mathfrak^{2b/\delta}_{K_x}}
\nc{\Ma}{\mathfrak^{2a/\delta}_{K_x}}
\nc{\bm}{\mf{N}^2}
\nc{\bn}{\mf{A}^2}
\nc{\al}{\alpha}
\nc{\barpsi}{\bar\psi}
\nc{\tldpsi}{\tilde\psi}
\nc{\tldphi}{\tilde\phi}
\nc{\dif}{\mathrm{d}}
\nc{\G}{\Gamma}
\nc{\g}{\gamma}
\nc{\bs}{\backslash}
\nc{\D}{\nabla}
\nc{\p}{\partial}
\nc{\DD}{\Delta^2}
\nc{\pp}{\partial^2} 
\nc{\de}{\delta}
\nc{\td}[2]{\trdeg{({#1}/{#2})}}
\nc{\dtd}[2]{\trdeg_{\delta}{({#1}/{#2})}}
\nc{\dspec}[1]{\spec_{\delta}{#1}}
\nc{\gens}[1]{\langle {#1} \rangle}        
\nc{\gen}[2]{ {#1} \langle {#2} \rangle } 
\nc{\form}{\Omega}
\nc{\set}[1]{\left\{ {#1} \right\}}
\nc{\mr}{\hat}
\nc{\pr}{\partial}
\nc{\bc}[3]{\cB^{#1}({#2},{#3})=B^{#1}_{#2}(C^{#2}_{#3}(t)+B^{#1}_{#3}(C^{#2}_{#3}(t))} 
\nc{\tuple}[2]{{#1},\ldots,{#2}} \nc{\ptu}[2]{{#1}:\ldots:{#2}}

\nc{\maps}[3]{{#1}\!:\!{#2}\rightarrow{#3}}
\nc{\map}[2]{{#1}\rightarrow {#2}} \nc{\res}[2]{{#1} |_{#2}}
\nc{\imbed}{\hookrightarrow}

\title{Quasi-conformal Rigidity of Negatively Curved Three Manifolds}
\author{ Yong Hou} 
\date{}
\begin{document}
\maketitle

\begin{abstract}
In this paper we study the rigidity of infinite volume $3$-manifolds with sectional 
curvature $-b^2\le K\le -1$ and finitely generated fundamental group.
In-particular, we generalize the 
Sullivan's 
quasi-conformal rigidity for finitely generated fundamental group with empty dissipative
set to negative variable curvature $3$-manifolds. We also generalize the 
rigidity of Hamenst\"{a}dt or more recently Besson-Courtois-Gallot, to 
$3$-manifolds with infinite volume 
and geometrically infinite fundamental group. 
\end{abstract}
\section{Introduction}
Let $\widetilde M$ be a simply connected complete Riemannian manifold with 
sectional curvature $-b^2\le K\le -1$. Let $\iso(\widetilde M)$ denote
the group of isometries of $\widetilde M$. Let $\G$ be 
a non-elementary, torsion-free, discrete subgroup of
$\iso(\widetilde M)$, and set $M:=\widetilde M/\G$.\par

First we recall some terminologies that is required for the statement of
the theorem. 
Let $S_\infty$ denote the boundary of $\widetilde M$. 
On $S_\infty$ one can define a metric in the following way. Let $v$ be a vector in the unit tangent bundle $S\widetilde M$. The geodesic
$v(t)$ defines two points on $S_\infty$ given by $v(\infty)$ 
and $v(-\infty)$. Let $\pi_t$ be the projection of $S_\infty\bs v(-\infty)$
along the geodesics which are asymptotic to $v(-\infty)$ to 
the horosphere which is tangent to $v(-\infty)$ and passing through $v(t)$.
Let $\dis_{v,t}$ be the distance on the horosphere induced by restriction of 
the Riemannian distance, $\dis$.
On $S_\infty\bs v(-\infty)\times S_\infty\bs v(-\infty)$ define
a function $\eta_v$ as $\eta_v(\xi,\zeta):=e^{-l_v(\xi,\zeta)}$ with
$l_v(\xi,\zeta):=\sup\{t| \dis_{v,t}(\pi_t(\xi),\pi_t(\zeta))\le 1\}$.
By our curvature assumption $-b^2\le K\le -1$, the function $\eta_v$ 
is a \emph{distance} on $S_\infty\bs v(-\infty)$, see \cite{Hem1}.\par

Every element of $\g\in\G$ has either exactly one or two fixed 
points in $S_\infty$, and $\g$ is called loxodromic if it has 
two fixed points \cite{BGS}. 
The group $\G$ is called
\emph{purely loxodromic} if all $\g\in\G$ are loxodromic. 
The limit set of $\G$ denoted by $\Lambda_\G$ is the unique minimal closed 
$\Gamma$-invariant subset of $S_\infty$ \cite{Gr2}. 
If $\G$ is purely loxodromic and
$\Lambda_\G=S_\infty$, then it can be either cocompact or 
$\widetilde M/\G$ is \emph{geometrically infinite}, hence $\G$ has infinite co-volume.
The \emph{convex hull} $CH_\G$ is the smallest convex set in $\widetilde M\cup S_\infty$
containing $\Lambda_\G$. The group $\G$ is called \emph{convex-cocompact}
if $CH_\G/\G$ is compact.\par
The critical exponent of $\G$ is the unique positive number $D_\G$
such that the Poincar\'{e} series of $\G$ given by 
$\sum_{\g\in\G}e^{-s\dis(x,\g x)}$ is divergent for $s<D_\G$ and
convergent for $s>D_\G$. If the Poincar\'{e} series diverges at
$s=D_\G$ then $\G$ is called \emph{divergent}.\par

Let $f:(X,\rho_X)\longrightarrow (Y,\rho_Y)$ be a embedding between two 
topological metric spaces. Then
$f$ is called \emph{quasi-conformal} embedding \cite{Y3} if there exists a constant 
$\kappa>0$ such that, for any $x\in X$ and $r>0$ there is $r_f(x,r)>0$ with 
\[f(X)\cap B'(f(x),r_f(x,r))\subset f(B(x,r))\subset B'(f(x),\kappa r_f(x,r)).\] 
where $B$ and $B'$ denotes a ball in $X$ and $Y$ respectively. 
When $f(X)=Y$ then $f$ is a quasi-conformal homeomorphism.\par

A torsion-free discrete subgroup $\G$ of $\iso(\widetilde M)$ is called
\emph{topologically tame} if $\widetilde M/\G$ is homeomorphic to the 
interior of a compact manifold-with-boundary. 

\begin{thm}\label{main-2}
Let $\G'\subset\Psl(2,\mathbb{C})$  be a topologically tame 
discrete group with $\Lambda_{\G'}=S^2$, and 
isomorphic $\chi:\G'\longrightarrow\G$ to a 
convex-cocompact discrete subgroup $\G$ of $\iso(\widetilde M)$ (here $\widetilde M$ is
$n$-dimensional). Let $f:S^2\longrightarrow S_\infty$ be a quasi-conformal
embedding which conjugate $\G'$ to $\G$, i.e. $f\circ\g=\chi(\g)\circ f$, for $\g\in\G'$. 
Then $D_\G\ge D_{\G'}$, and equality
if and only if $\mathbb{H}^3$ embeds isometrically into $\widetilde M$ and the
action of $\G$ stabilizes the image.
\end{thm}
To state our next theorem we need to introduce one additional terminology.
We take $\widetilde M$ to be a $3$-manifold in the following.\par
Let $\mf{M}^\lambda_{\eta_v}$ denote the $\lambda$-dimensional hausdorff measure
on $(S_\infty\bs v(-\infty),\eta_v)$. We say $\G$ is 
\emph{hausdorff-conservative} if there exists a constant $\alpha(v)>0$
such that 
$\alpha^{-1}r^{D_\G}\le\mf{M}^{D_\G}_{\eta_v}(B(\xi,r)\cap\Lambda_\G)\le\alpha r^{D_\G}$ 
for any ball $B(\xi,r)$ of
radius $r$ about $\xi\in\Lambda_\G$ in $(S_\infty\bs v(-\infty),\eta_v)$. 
From this definition, we note that if $\G$ is a finitely generated 
torsion-free 
discrete subgroup of $\Psl(2,\mathbb{C})$ with $D_\G=2$, then 
hausdorff-conservative implies \emph{conservative} 
(classical definition, $\S 5$). Conversely, if $\G$ 
is a topologically tame, \emph{conservative}, discrete subgroup of 
$\Psl(2,\mathbb{C})$, then $\G$ is hausdorff-conservative, 
see Proposition \ref{hd-conservative}. We believe all finitely generated 
conservative discrete subgroup of $\Psl(2,\mathbb{C})$ are 
hausdorff-conservative, see Remark \ref{conser-hd}.
For a convex-cocompact $\widetilde M/\G$ with $-b^2\le K\le -1$, 
it follows from \cite{Coor}, $\G$ is hausdorff-conservative. Now we are
ready to state the theorem which generalizes Sullivan's quasi-conformal 
rigidity theorem.
\begin{thm}[Main]\label{main-1-2}
Let $\G$ be a topologically tame, purely loxodromic discrete subgroup
of $\iso(\widetilde M)$ with $\Lambda_\G=S_\infty$. Let $\G'$ be a topologically tame  discrete subgroup
of $\Psl(2,\mathbb{C})$. Suppose $f:S_\infty\longrightarrow S^2$ is
a quasi-conformal homeomorphism conjugate $\G$ to $\G'$. 
Then $D_\G\ge D_{\G'}$, and $\G=\g\G'\g^{-1}$ with
$\g\in\Psl(2,\mathbb{C})$ if and only if $D_\G=D_{\G'}$ and 
$\G$ is hausdorff-conservative.
\end{thm}
\begin{cor}\label{cor-1}
Let $M=\widetilde M/\G$ be a complete topologically tame $3$-manifold 
with $-b^2\le K\le -1$, $\G$ purely loxodromic, 
and $\Lambda_\G=S_\infty$. 
Let $h:M\longrightarrow N$ be a quasi-isometric homeomorphism 
to a hyperbolic manifold $N$. Then $M$ is isometric to $N$ 
if and only if $D_\G=2$ and $\G$ is hausdorff-conservative.
\end{cor}
Let us point out that Theorem \ref{main-1-2}, generalizes known
rigidity theorems in two directions for three dimensional manifolds.\par
First assume $M$ is hyperbolic ($b=1$) but not necessarily 
geometrically finite. Since $M$ is topologically tame and $\Lambda_\G=S^2$ 
we have $D_\G=2$ by analytical tameness (see Proposition \ref{lem-1}). 
Hence by Theorem \ref{main-1-2}, $M$ is
quasi-conformal stable. This is a case of the Sullivan rigidity 
theorem for topologically tame $\G$ with empty dissipative set.
Next let us assume $M$ is compact with $-b^2\le K\le -1$.
Then the critical exponent $D_\G$ is equal to $h_M$ the topological entropy of $M$, 
and by \cite{Gabai}, any homotopy equivalence between $M$ and a compact
hyperbolic $3$-manifold is induced by a homeomorphism.
Therefore it follows from Corollary \ref{cor-1} we have: $M$ is isometric to a compact hyperbolic $3$-manifold 
if and only if they are homotopically equivalent and $h_M=2$.
This is the Hamenst\"{a}dt's rigidity or more recently Besson-Courtois-Gallot
theorem for $3$-manifolds.\par
Note that it also follows from Theorem \ref{main-1-2}, the 
quasi-conformal version 
of the Hamenst\"{a}dt's theorem for compact $3$-manifold $M$  can be stated as:
\begin{cor}\label{main-3-3}
Let $\G$ be a cocompact discrete subgroup
of $\iso(\widetilde M)$. 
Let $\G'\subset\Psl(2,\mathbb{C})$ be a discrete group. 
Suppose $f:S_\infty\longrightarrow S^2$ is a
quasi-conformal homeomorphism conjugate $\G$ to $\G'$, 
Then $D_\G\ge D_{\G'}$, and equality if and only if $\widetilde M/\G$ 
is isometric to $\mathbb{H}^3/\G'$.
\end{cor}
The proves of these theorems relies on our next result,
\begin{thm}\label{Main4}
Let $M=\tilde{M}/\G$ be a topologically tame 
$3$-manifold with $-b^2\leq\cK\leq -1$. Suppose that $\G$ is purely loxodromic and that $\Lambda(\G)=S_\infty$. 
Then $2\le D$ and $\G$ is harmonically ergodic. If $D=2$ then $\G$ 
is also divergent. 
\end{thm} 
In section $2$, we state some of the topological properties of negatively pinched 
$3$-manifolds. In particular, we define \emph{geometrically infinite ends} for 
negatively pinched 
$3$-manifolds, and then state our theorem which describe the geometrical properties of
this type of end, it is a crucial step in the proof of Theorem \ref{Main4}.
Section $3$ discusses measures on $S_\infty$ and the ergodicity of $\G$ 
with respect to these measures. In section $4$, we give proofs of part I of the
theorems. And section $5$ is used to complete the proofs.\\
\centerline{ ACKNOWLEDGMENTS}\par 
I am very grateful to Peter Shalen for his unwavering support and encouragement. 
I wish to thank Marc Culler, Dick Canary for their interest. 
\section{Topological Ends} 
Every isometry of $\tilde{M}$ can be extend to a Lipschitz map on 
 $S_\infty:=\partial\tilde{M}$ \cite{Gr2}. For a torsion-free $\G$, every element 
$\g\in\G$ is one of the following types: $(1)$ \emph{parabolic} if it has 
exactly one fixed point in $\tilde{M}\cup S_\infty$ which lies in
$S_\infty$; $(2)$ \emph{loxodromic} if it has exactly two distinct fixed 
points in $\tilde{M}\cup S_\infty$, both lying in $S_\infty$. \par
Denote by $\Lambda (\Gamma)\subset \p\tilde M$ the limit set of $\Gamma$, 
which is the unique minimal closed $\Gamma -$invariant subset of $S_\infty$. 
Most of the important properties of the limit set in the constant curvature space continue 
to hold in the variable curvature space \cite{EO}. In particular: (i) $\Lambda(\G)=\overline{\G x}\cap S_\infty$; (ii) 
$\Lambda(\G)$is the closure of the set of fixed points of loxodromic elements of 
$\G$; and (iii) $\Lambda(\G)$ is a perfect subset of $\G$. The set 
$\Omega(\G):=S_\infty\backslash\Lambda(\G)$ is the region of discontinuity. 
The action of $\G$ on $\tilde{M}\cup\Omega(\G)$ is proper and discontinuous, 
see \cite{EO}. The manifold $M_\G:=\tilde{M}\cup\Omega(\G)/\G$ with possibly 
nonempty boundary is traditionally called the Kleinian manifold. 
We also let $\Lambda_c(\G)$ denote the conical limit set of $\G$, 
i.e. $\xi\in\Lambda_c(\G)$ if for some $x\in\tilde M$ (and hence for every $x$) there exist a sequence $(\g_n)$ of elements in 
$\G$, a sequence $(t_n)$ of real numbers, and a real number $C>0$, such that $\g_n x\longrightarrow\xi$ and 
$\dis(c^\xi_x(t_n),\g_n x)<C$ where $c^\xi_x$ is the geodesic ray connecting
$x$ and $\xi$. Equivalently, a point belongs to $\Lambda_c(\G)$ if it belongs to infinitely many shadows cast by 
balls of some fixed radius centered at points of a fixed orbit of $\Gamma$. Note that $\Lambda_c(\G)$ is a $\G-$invariant 
subset of $\Lambda(\G)$, hence a dense subset.
\begin{prop}[Margulis Lemma]\label{Margulis}
There exists a number $\epsilon_b$ which only depend 
on the pinching constant $b$ of $M$, such that the 
group $\G_{\epsilon}$ generated by elements in $\Gamma$ of length at 
most $\epsilon_b$ with respect to a fixed 
point in M is almost nilpotent of rank at most $2$. Then the number, 
$2\epsilon_b$ is called the \emph{Margulis constant}.
\end{prop}
Note that, if $M$ is orientable and 
$\G$ is torsion-free, then Margulis Lemma implies $\G_{\epsilon_b}$ is abelian.\par
Let $\epsilon\le\epsilon_b$ be given. Then $M$ may be written as the 
union of a {\it thin part}  $M_{[0,\epsilon)}$ consisting of all points at 
which there is based a homotopically nontrivial loop of length $\le\epsilon$ 
and a {\it thick part}  $M_{[\epsilon,\infty)}=\overline{ M-M_{[0,\epsilon)}}$. 
Note that $M_{[\epsilon,\infty)}$ is compact if $M$ is of finite volume. 
Also the thin part of $M$ is completely classified by the next proposition.
\begin{prop}\label{cusps}
Each connected component of $M_{[0,\epsilon)}$ is diffeomorphic to one 
of the following :
\begin{description}
\item[parabolic rank-1 cusp]: $S^1\times{\mathbb R}\times [0,\infty)$.
\item[parabolic rank-2 cusp]: $T^2\times [0,\infty)$.
\item[solid torus about the axis of a loxodromic $\gamma$]: $D^2\times S^1$. 
\end{description}
\end{prop}
 For simplicity we restrict to the case where $M$ has no cusps. It follows 
from the existence of a  compact core $C(M)$ for $M$ \cite{FM}, that $M$ has 
only finitely many ends \cite{BP}. In fact, each component of 
$\partial C(M)$ is the boundary of a neighborhood of an end of $M$, and 
this gives a bijective correspondence between ends of $M$ and components 
of $\partial C(M)$.\par
We define the simplicial ruled surfaces as 
follows. Let $S$ be a surface of positive genus and let $T_P$ be a 
triangulation defined with respect to  a finite collection $P$ of points 
of $S$. This means that  $T_P$ is a maximal collection of nonisotopic 
essential arcs with end points in $P$; these arcs are the {\it edges} of 
the triangulation, and the components of the complement in $S$ of the 
union of the edges are the {\it faces}. Let $f \text{ : } S\longrightarrow 
M$ be a map which takes edges to geodesic arcs and faces to nondegenerate 
geodesic ruled triangles in $M$.  The map $f$ induces a singular metric on 
$S$. If  the total angle about each vertex of $S$ with respect to this 
metric is at least $2\pi$, then the pair $(S,f)$ is called a 
\emph{simplicial ruled surface}. It follows from the definition of the 
induced metric on $S$ that $f$ preserves lengths of  paths and is therefore 
distance non-increasing. Any geodesic ruled triangle in $M$ has Gaussian 
curvature at most $-a^2$. This means that each $2$-simplex of $S$ inherits 
a Riemannian metric of curvature at most $-a^2$. Since we have required 
the the total angle at each vertex to be at least $2\pi$, by Gauss-Bonnet 
theorem the curvature of $S$ is negative in the induced metric.
\begin{defn}\label{Geo-inf-def}
An end $E$ is said to be a \emph{geometrically infinite} if there exists a divergent
sequence of geodesics, i.e: there exists a sequence of closed geodesics
${\alpha_k}\subset M^{\circ}_\epsilon$, such that for any neighborhood $U$ of
$E$, there exists some positive integer $N$ such that $\alpha_k\subset U$ for all $k>N$. 
If in addition for some surface $S_E$ we have that $U$ is homeomorphic 
to $S_E \times [0,\infty)$, and there exists a sequence
of simplicial ruled surfaces : $S_E\stackrel{f_l}{\longrightarrow}U$ such that $f_l(S_E)$ is homotopic to $S_E\times {0}$ in $U$ and leaves every
compact subset of $M$, then $E$ is said to be \emph{simply degenerate}. The sequence
$(S_E\stackrel{f_l}{\longrightarrow}U)$ is called an \emph{exiting sequence}.
A end which is not geometrically infinite will be called 
\emph{geometrically finite}. 
\end{defn} 
\begin{thm}[Hou]\label{a16} 
Let $M=\tilde{M}/\G$ be a topologically tame negatively pinched $3$-manifold with $\G$ purely
loxodromic. Then all geometrically infinite ends of $M$ are simply degenerate.
And if $\Lambda(\G)=S_\infty$, then there 
are no nonconstant positive superharmonic
functions, or nonconstant subharmonic functions bounded above, on $M$.
\end{thm}
\section{$\G$-action}
In this section we will study the action of $\G$ on $S_\infty$ and prove 
ergodicity of $\G$ for topologically tame $3$-manifolds with $\Lambda(\G)=S_\infty$. 
We will prove that for such a manifold, the Green series is divergent, and 
that the Poincar\'{e} series is also divergent if $D=2$. 
Theorem \ref{Main4} will also be 
proved in this section.\par
In some situations we will take the dimension of $M$ to be $3$, otherwise 
we will assume $M$ is $n$-dimensional in general.\par

Set the following notations throughout the paper.
Let $\G'\subset\Psl(2,\mathbb{C})$ be a discrete torsion-free subgroup.
Denote $S^2:=\pr\mathbb{H}^3$, and $S_\infty:=\pr\widetilde M$.

There are many equivalent ways of equipping $S_\infty$ with a metric which is 
compatible with $\G$-action. 
Fix a point $x\in\widetilde M$. Let $\xi,\zeta$ in $S_\infty$ be given. Set
$c^\zeta_y(t)$ as the geodesic ray connecting $y$ and $\zeta$.\par
In \cite{Gr2}, Gromov defined a metric on $S_\infty$ as follows. For $y,z\in\widetilde M$, let 
us consider arbitrary
 continuous curve $c(t)$ in $\widetilde M$ with initial point and end point denoted by $c(t_0)=y$ and
$c(t_1)=z$ respectively. Define a nonnegative real-valued function $\cG_x$ on 
$\widetilde M\times\widetilde M$ by 
\[ \cG_x(y,z):=\inf_{\text{all }c}\left(\int_{[t_0,t_1]}e^{-\dis(x,c(t))}\dif t\right).\]
In particular, Gromov showed 
the function $\cG_x$ extends continuously to $S_\infty\times S_\infty$.
Every element of $\G$ extends to $S_\infty$ as a Lipschitz map with respect
to $\cG_x$.\par 
In \cite{Ka} the following metrics are shown to be equivalent to the Gromov's metric. 
\begin{description}
\item[$K_x$ metric]:
Let $B_\zeta$ denote the Busemann function based at $x_0$.    
Set $B_\zeta(x,y)$\\$=B_\zeta(x)-B_\zeta(y)$, for $x,y \in\widetilde M$, 
the function $B_\zeta(x,y)$ is called the Busemann cocycle.
Define $\beta_x:$ $S_\infty\times S_\infty
\longrightarrow\mathbb{R}$ by $\beta_x(\xi,\zeta):=B_\xi(x,y)+B_\zeta(x,y)$
where $y$ is a point on the geodesic connecting $\xi$ and $\zeta$. 
The $K_x$ metric is then defined by
\[ K_x(\xi,\zeta):=e^{-\frac{1}{2}\beta_x(\xi,\zeta)}. \]
\item[$L_x$ metric:]
Let $\alpha_x(\xi,\zeta)$ denote the distance between $x$ and the geodesic connecting
$\xi$ and $\zeta$. The function $L_x:$ $S_\infty\times S_\infty\longrightarrow
\mathbb{R}$ is then defined by
\[ L_x(\xi,\zeta):=e^{-\alpha_x(\xi,\zeta)}. \]
\item[$d_x$ metric:]
Define a function $l_x:$ $S_\infty\times S_\infty\longrightarrow\mathbb{R}$ by
$l_x(\xi,\zeta):=\sup\{\tau|$ $\dis(c^x_\xi(\tau),c^x_\zeta(\tau))=1\}$.
 Geometrically, a neighborhood about $\xi$ in $S_\infty$ with respect to 
the topology induced by $l_x$ is the shadow cast by the intersection of 1-ball about $c^\xi_x(\tau)$
and $\tau$-sphere about $x$. The $d_x$ metric is then defined by
\[ d_x(\xi,\zeta):=e^{-l_x(\xi,\zeta)}. \]  
\end{description}
It was originally observed for symmetric spaces by Mostow \cite{Mos1} that the boundary map is 
quasi-conformal. This property continue 
to hold in negatively curved spaces, see \cite{Hem1} and \cite{Pansu}. Here we give a 
proof of this fact with respect to the above metrics. 
\begin{prop}\label{qc-1}
Let $h$ be a quasiisometry between two negatively pinched curved spaces. 
The boundary extension map $\bar h$ is quasi-conformal on the boundary with 
respect to $d_x,L_x,K_x,\eta_v$-metrics.
\end{prop}
\begin{proof}
For the proof of $\eta_v$-metric See Proposition $3.1$ in \cite{Hem1}. 
Fix $x\in\tilde M$. Let us take $d_x$-metric. Set $\lambda\ge L$. Denote by
$S(x;y,R)$ the shadow cased from $x$ of the metric sphere $S(y,R)$ with center 
located at $y$ and radius $R$, 
i.e. $S(x;y,R)=\{\xi\in S_\infty|c^\xi_x\cap S(y,R)\not=\emptyset\}$. 
Let $B(\xi,r)$ be
a ball of radius $r$ in $S_\infty$. Using triangle comparison we can show
there exists a constant $\alpha_b\ge 1$ depends on pinching constant $b$ such 
that 
\[S(x;c^\xi_x(t_r),\lambda)\subset B(\xi,r)
\subset S(x;c^\xi_x(t_r),\alpha_b\lambda)\]
for some $t_r>0$ which depends only on $r$. The images 
$\bar\phi(S(x;c^\xi_x(t_r),\lambda))$ and $\bar\phi(S(x;c^\xi_x(t_r),\alpha_b\lambda))$
are quasi-spheres, i.e. there exists a constant $\beta_\phi>0$ depends on $\phi$
such that $S(\bar\phi(x); \bar\phi(c^\xi_x(t_r)),\beta^{-1}_\phi\lambda)
\subset\bar\phi(S(x;c^\xi_x(t_r),\lambda))$ and 
$\bar\phi(S(x;c^\xi_x(t_r),\alpha_b\lambda)\subset S(\bar\phi(x);\bar\phi(
c^\xi_x(t_r)),\beta_\phi\alpha_b\lambda)$. On the other hand, by estimates
in \cite{C1} there exists positive numbers $A_1(\beta_\phi,\lambda)$
and $A_2(\alpha_b,\beta_\phi,\lambda)$ such that
\[B(\bar\phi(\xi),A_1e^{-R})\subset S(\bar\phi(x); \bar\phi(c^\xi_x(t_r)),\beta^{-1}_\phi\lambda),\]
\[S(\bar\phi(x);\bar\phi(c^\xi_x(t_r)),\beta_\phi\alpha_b\lambda)\subset
B(\bar\phi(\xi),A_2e^{-R})\] where $R=\dis(\bar\phi(x),\bar\phi(c^\xi_x(t_r))$.
Hence the result follows by setting $r_\phi(\xi,r)=A_1e^{-R}$ and
$\kappa=A_2/A_1$.
\end{proof}

\begin{prop}\label{prop-1}
Let $f:\pr\widetilde N\longrightarrow S_\infty$ be a 
embedding conjugate $\G_1$ to $\G_2$ 
under isomorphism 
$\chi:\G_1\longrightarrow\G_2$ 
($f \circ\g=\chi(\g)\circ f$).
Then $f(\Lambda_{\G_1})=\Lambda_{\G_2}$.
\end{prop}
\begin{proof}
Let $\g\in\G_1$.
Since 
$\g f^{-1}(\Lambda_{\G_2})=f^{-1}(\chi(\g)\Lambda_{\G_2})$, 
and by $\G_2$-invariance of $\Lambda_{\G_2}$, we have 
$f^{-1}(\Lambda_{\G_2})$ is $\G_1$-invariant closed set. Note that
$f^{-1}(\Lambda_2)$ is nonempty, since fixed 
points of elements of $\G_1$ are also fixed points of elements of $\G_2$,
hence $f(\Lambda_{\G_1})\subseteq\Lambda_{\G_2}$. Similarly 
we also have $f(\Lambda_{\G_1})\supseteq\Lambda_{\G_2}$, and result follows. 
\end{proof}
\begin{prop}\label{lem-1}
Let $\G$ be a topologically tame, torsion-free, discrete subgroup of
$\iso(\widetilde M)$ with 
$\Lambda_\G=S_\infty$. Let $\G'$ be a topologically tame, discrete subgroup
of $\Psl(2,\mathbb{C})$. Suppose $f:S_\infty\longrightarrow S^2$
is a homeomorphism conjugate $\G$ to $\G'$. 
Then $D_{\G'}=2$ and $\G'$ is divergent. 
\end{prop}
\begin{proof}
By Proposition \ref{prop-1} the hyperbolic manifold $N=\mathbb{H}^3/\G'$ is
topologically tame and $\Lambda_{\G'}=S^2$. It follows from
analytical tameness and Theorem $9.1$ of \cite{Ca}, there exists no
non-trivial positive superharmonic function on $N$ with respect to the hyperbolic Laplacian
$\Delta$. Let $P(y,\xi)$ denote the Poisson kernel on $\mathbb{H}^3$. 
The $D_{\G'}$-dimensional conformal measure 
(Patterson-Sullivan measure, see end of $\S 3$) $\sigma_y$ has Radon-Nikodym 
derivative of
$P(y,\xi)^{D_{\G'}}$, i.e. $\frac{\dif\sigma_y}{\dif\g^*\sigma_y}(\xi)=P(\g^{-1}y,\xi)^{D_{\G'}}$.
The $\G'$-invariant function $h(y):=\sigma_y(S^2)$ satisfies 
$\Delta h=D_{\G'}(D_{\G'}-2)h$, which implies $h$ is non-trivial superharmonic if 
$D_{\G'}\not= 2$. And it follows that $\G'$ must also be divergent.
\end{proof}
Let $C$ be a subset of $S_\infty$. Let $\lambda$-dimensional Hausdorff measure 
of $C$ on the metric space$(S_\infty,\rho_x)$ be denoted by 
${\mf M}^\lambda_{\rho_x}(C)$. 
Observe that for any $x\in\widetilde M$ and any $\gamma\in\G$, 
we have $\gamma^*{\mathfrak M}^\lambda_{K_x}={\mathfrak M}^\lambda_{K_{\gamma^{-1}x}}$; 
this follows from the straightforward identity.\par
A family of finite Borel measures $[\nu_y]_{y\in\widetilde M}$.
will be called a 
$\lambda$-{\it conformal density under the action of} $\G$ 
if for every $x\in\widetilde M$ and every $\gamma\in\G$ we have 
$\gamma^*\nu_y=\nu_{\gamma^* y}$, and 
the Radon-Nikodym derivative $\frac{\dif \nu_y}{\dif\gamma^*\nu_y}(\zeta)$  
at any point $\zeta\in S_\infty$ is equal to 
$e^{-\lambda B_\zeta(\gamma^{-1}y,y)}$. (This is to be interpreted as 
being vacuously true if, for example, the measures in the family are all 
identically zero). Although there can not be any $\G$-invariant non-trivial finite
Borel measure on $\Lambda_\G$ for non-elementary $\G$, we can always 
define a $\G$-invariant non-trivial locally finite measure $\Pi_{\nu_x}$ on 
$\Lambda_\G\times\Lambda_\G$ by setting 
$\dif\Pi_{\nu_x}(\xi,\zeta)=e^{\lambda\beta_x(\xi,\zeta)}\dif\nu_x(\xi)\dif\nu_x(\zeta)$.
The measure $\Pi_{\nu_x}$ corresponds to the Bowen-Margulis measure on the 
unit tangent bundle $S\widetilde M$ see \cite{Ka}.\par

Let us recall a fundamental fact about conformal density, which was 
originally proved by Sullivan for $\G\subset\SO(n,1)$ and generalized to 
the pinched negatively curved spaces in \cite{Y2}. It relates the divergence 
of $\Gamma$ at the critical exponent $D_\G$ with ergodicity of the 
$D_\G$-conformal density under the action of $\Gamma$.\par
We will say that two Borel measures on $S_\infty$ are in the same 
$\G$-class if the Radon-Nikodym derivative of $\g^*\nu_1$ with respect 
to $\nu_1$ is equal to the Radon-Nikodym derivative of $\g^*\nu_2$ with 
respect to $\nu_2$.
\begin{prop}[see \cite{Y2}]\label{conf}
Let $\Gamma$ be a nonelementary, discrete, torsion-free and divergent 
at $D_\G$. Suppose $[\nu]$ is a $D_\G$-conformal density under the action of 
$\Gamma$, then $\Gamma$ act ergodically on $\Lambda_\G$ and 
$\Lambda_\G\times\Lambda_\G$ with respect to $[\nu_x]$ and 
$[\Pi_{\nu_x}]$ respectively.
\end{prop}
\begin{prop}[see \cite{Nil}]\label{prop-unique}
Let $\Gamma$ be nonelementary and discrete. Suppose that $\Gamma$ acts ergodically
on $S_\infty$ with respect to a measure $\nu$ defined on $S_\infty$. 
Then every measure of $S_\infty$ in the same measure class as $\nu$ is a constant
multiple of $\nu$.
\end{prop}
\begin{prop}\label{a34}
Let $\Gamma$ be a non-elementary discrete subgroup of the isometry group
of $\tilde{M}$. If $[\nu_y]^D_{y\in\tilde{M}}$ is a non-trivial 
$\G$-invariant $D$-conformal density, then $D\not= 0$.
\end{prop}
\begin{proof}
Suppose $D=0$. Then $\nu_y$ is a $\G$-invariant non-trivial finite Borel 
measure. Since $\G$ is non-elementary, there exists a loxodromic element 
$\g$ in $\G$. Let $\xi,\zeta\in S_\infty$ be the two distinct fixed points 
of $\g$. Let $<\g>$ be the group generated by $\g$. Then $\nu_y$ is clearly 
$<\g>$-invariant. But $\g$ is loxodromic, so we must have 
$\text{supp}(\nu_y)\subset\{\xi,\zeta\}$. Then, by the fact that
$\Lambda(\G)$ is infinite, we have $\nu_y$ is an infinite measure, which is 
a contradiction.
\end{proof}
\begin{prop}\label{a26} 
Let $\Gamma$ be a discrete subgroup of $\iso(\widetilde M)$. Suppose 
${\mathfrak M}^\lambda_{K_x}$ is a finite measure. Then 
${\mathfrak M}^\lambda_{K_x}$ 
is a $\lambda$-conformal density under the action of $\Gamma$.
\end{prop}
There is a canonical way of constructing $D_\G$-dimensional conformal
density which is due to Patterson-Sullivan as follows; By applying
a adjusting function we can always assume the Poincar\'{e} series
diverges at $D_\G$. The measures
\[\mu_{x,s}:=\frac{\sum_{\g\in\G}e^{-s\dis(z,\g x)}\delta_{\g x}}
{\sum_{\g\in\G}e^{-s\dis(z,\g z)}}\quad;s>D_\G\]
converges weakly to a limiting measure $\mu_x$ as $s_n\to D_\G$ through
a subsequence. It is trivial to see that $\mu_x$ is supported on
$\Lambda_\G$. The measure $[\mu_x]$ is called Patterson-Sullivan measure 
which is $D_\G$-conformal under $\G$ see \cite{Pat}, \cite{Sul2}.\par
From now on for $\G'\subset\Psl(2,\mathbb{C})$, we will denote the 
Patterson-Sullivan measure on $\Lambda_{\G'}$ by $\sigma_y$.\par
Let $\Lambda^c_\G$ denote the set of conical limit points in $\Lambda_\G$. 
Recall a point $\xi\in\Lambda_\G$ is in $\Lambda^c_\G$ if and only if 
there exists $\{\g_n\}\subset\G$ such that $\dis(\g_n c^\xi_x(t_n),x)<c$ for some
$c>0$ and sequence of $t_n$. Obviously $\Lambda^c_\G$ is $\G$-invariant, and
non-empty (a loxodromic fixed points are in $\Lambda^c_\G$), hence 
it is a dense $\G$-invariant subset of $\Lambda_\G$.
A equivalent definition for the conical limit
point $\xi$ is that it must be contained in infinitely many shadows 
$S(x;\g_n x,c)$. Hence 
$\Lambda^c_\G=\cup_{\lambda>0}\cap_{m\ge 1}\cup_{n>m}S(x,\g_n x,\lambda)$.  
It is a easy fact from the construction
of $\mu_x$, no points in $\Lambda^c_\G$ can be a atom for $\mu_x$,
and if $\supp(\mu_x)\subseteq\Lambda^c_\G$ then $\G$ is divergent.
In fact it is a deep result of Sullivan that $\G$ is divergent if and only
if $\supp(\mu_x)\subseteq\Lambda^c_\G$. 
\begin{lem}[see \cite{Coor}; Sullivan's Shadow Lemma]\label{shadow}
Let $\mu_x$ be a $D_\G$-conformal density with respect to $\G$, which is
not a single atom. Then there exists constants $\al>0$ and $\lambda_o\ge 0$,
such that,
\[\al^{-1}e^{-D_\G\dis(x,\g^{-1}x)}\le\mu_x(S(x;\g x,\lambda))
\le\al e^{-D_\G\dis(x,\g^{-1} x)+2D\lambda},\]
for all $\g\in\G$ and $\lambda\ge\lambda_o$.
\end{lem}

\begin{prop}\label{no-atom}
Let $\G\subset\iso(\widetilde M)$ be a discrete subgroup. 
Suppose either $\Lambda_\G=\Lambda^c_\G$ or $\Lambda_\G=S_\infty$ and 
$\G$ is divergent. 
Then $\mu_x$ is positive on all non-empty relative open subsets 
of $\Lambda_\G$.
\end{prop}
\begin{proof}
Suppose $\Lambda_\G=S_\infty$. It suffices to show $\mu_x$ is positive for any non-empty open ball 
$B(\xi,r)$ with respect to the $d_x$-metric. Fix $\lambda>\lambda_o$. 
Let $\zeta\in\Lambda^c_\G\cap B(\xi,r)$ (note that $\Lambda^c_\G$ is dense in
$\Lambda_\G$ so the intersection is nonempty). Then we can choose $\gamma\in\G$ such that 
$S(x;\gamma x,\lambda)\subset B(\xi,r)$. 
By assumption $\G$ is divergent, we have $\supp(\mu_x)\subseteq\Lambda^c_\G$. Since
no points of $\Lambda^c_\G$ can be a atom for conformal density, 
the result follows from Lemma \ref{shadow}. Same argument works if
$\Lambda_\G=\Lambda^c_\G$
\end{proof}

\textnormal{Let us define a function $\Theta :\tilde{M}\times\tilde{M}\times S_\infty
\longrightarrow\mathbb{R}^+$ by $\Theta(x,y,\xi):=\exp(-B_\xi(x,y))$.}\par
{\bf Harmonic Density}\\
\textnormal{Let $\lambda_1$ and $\tilde{\lambda}_1$ denote the first of the spectrum of $\Delta$ on $M=\tilde{M}/\G$,
and of $\tilde{\Delta}$ on $\tilde{M}$, respectively. Recall that for a noncompact
open manifold, the first of the spectrum is defined as 
\[ \lambda_1:=\inf_{f\in C^\infty_o,f\not=0}\left(\frac{\int |\nabla f|^2}
{\int f^2}\right) \] 
where $C^\infty_o$ is the space of smooth functions on $M$ with compact 
support. Note that we always have} 
$\lambda_1\le\tilde{\lambda}_1$.\par
\textnormal{The $\lambda_1$-harmonic functions has been studied by Ancona 
in \cite{Ancona} and \cite{Ancona2}.} 
\begin{prop}[Ancona]\label{a35}
For each $s<\lambda_1$, the elliptic operator $\tilde{\Delta}+s I$
has a Green function $G_s(x,y)$, and there exists a function
$f :\mathbb{R}^+\longrightarrow\mathbb{R}^+$ such that $\sum_{\g\in\G}\hat{G}_s
(x,\g y)$ converges for $s<\lambda_1$ and diverges for $s\geq\lambda_1$,
where $\hat{G}_s(x,\g y):=\exp(f(\dis(y,\g y)))G_s(x,\g y)$. Furthermore,
$\fP_s(x,y,\zeta):=\lim_{z\to\zeta}\frac{G_s(x,z)}{G_s(y,z)}$ defines the Poisson kernel
of $\tilde{\Delta}+s I$ at $\zeta\in S_\infty$.
\end{prop}
Similarly to the construction of $\mu_x$ from $Z_\G$, one can also construct a family of
Borel measures from $\sum_{\g\in\G}\hat{G}_s(x,\g y)$.
\begin{prop}\label{a36}
Let $x$ be any point of $\tilde M$. There exists a family of Borel 
measures $[\omega^1_y]_{y\in\tilde{M}}$ on $S_\infty$ such that (i) for all $x,y\in\tilde M$,
 Radon-Nikodym derivative ${\mathrm d}\omega_y^1/{\mathrm d}\omega_x^1$ at any point 
$\zeta\in S_\infty$ is equal to $\fP_{\lambda_1}(x,y,\zeta)$
and (ii) $\omega^1_x$ is of mass $1$.
\end{prop}
\indent\textnormal{Let us denote the harmonic density of} $\tilde{\Delta}$ \textnormal{by} 
$[\omega_y]_{y\in\tilde M}$ 
\textnormal{with $\omega_x$ normalized of mass $1$. By definition this means that every 
harmonic function $f$ on $M$ with boundary values $f_\infty$ is given by 
\[f(x)=\int_{S_\infty}f_\infty(\xi){\mathrm d}\omega_x(\xi).\]
 The existence and uniqueness of harmonic density
follows from the solvability of the Dirichlet problem on 
$\tilde{M}\cup S_\infty$
see \cite{AnSc} and the Riesz Representation Theorem. The Radon-Nikodym 
derivative of 
$\omega_x$ at $\xi\in S_\infty$ is given by the Poisson kernel 
$\fP(x,y,\xi)$ of $\tilde{\Delta}$, i.e. 
$\frac{{\dif}\omega_y}{{\dif}\omega_x}(\xi)=\fP(y,x,\xi)$. For any  $\G$-invariant subset
$C\subset S_\infty$, the function $h_C$ on 
$\tilde{M}$ defined by $h_C(y):=\int_{S_\infty}\chi_C\fP(y,x,\xi){\dif}\omega_x(\xi)$ 
is $\G$-invariant, hence defines a harmonic function on} $M$.
\begin{prop}\label{a37}
Let $M=\tilde M/\G$ be a negatively pinched topologically tame $3$-manifold  
with $\Lambda(\G)=S_\infty$. Then $\G$ is ergodic with respect to harmonic  
density $[\omega_y]_{y\in\tilde M}$.
\end{prop}
\begin{proof}
Suppose not, and let $C\subset S_\infty$ be a $\G$-invariant subset with $\omega_x(C)>0$
and $\omega_x(C^c)>0$. By Fatou's conical convergence theorem, we have
$\chi_C(\xi)=\lim_{t\to\infty}h_C(c^\xi_y(t))$ for $\xi\in S_\infty$. Hence, $h_C$
defines a positive nonconstant $\G$-invariant harmonic function, which 
contradicts Theorem \ref{a16}. Therefore
$\G$ must be ergodic.
\end{proof}
\begin{prop}\label{a38}
Let $M$ be noncompact and satisfy the hypothesis of Proposition \ref{a37}. Then $\omega^1_x=\omega_x$.
\end{prop}
\begin{proof}
Let us note that $\lambda_1=0$. This follows from the fact that for a non-compact, 
complete Riemannian manifold $M$, if $\lambda_1(M)>0$ then there exists a
positive Green's function $G$ on $M$. If such a $G$ exists, then $1-\exp(-G)$ defines a positive
superharmonic function, which is a contradiction to Theorem \ref{a16}. Hence we 
must have $\lambda_1=0$. Therefore $\fP_{\lambda_1}=\fP$, i.e $\frac{\dif\omega^1_y}{\dif\omega^1_x}
=\frac{\dif\omega_y}{\dif\omega_x}$. Hence, by Proposition \ref{a37} and uniqueness,
we have the desired result.
\end{proof}
{\bf Superharmonic Functions}\\
\textnormal{Let $\xi\in S_\infty$ be given. Let $\cE$ be a continuous unit 
vector field on $\tilde{M}$ with $\cE(x)=\Phi_x(\xi)$. Then by using the 
first length variation formula, one can 
show that $B_\xi$ is $C^1$ and that $-\grad B_\xi=\cE$. In fact, the Busemann function 
is $C^{2,\alpha}$, see \cite{HI}.}\par 
\textnormal{Let $[\mu_y]^D$ denote $D$-conformal density. Let us define a nonnegative function
$u$ on $\tilde{M}$ by} 
\[  u(y):=\int_{S_\infty} \Theta^D(y,x,\xi) {\mathrm d}\mu_x (\xi). \]
\begin{prop}\label{a39}
The function $u$ is a $\Gamma$-invariant and positive. It is superharmonic if 
$D\leq (n-1)a$, and subharmonic if 
$(n-1)a\le D\leq (n-1)b$.
\end{prop} 
\begin{proof}
We can write $u(y)$ as $\mu_y(S_\infty)$. Since
$u(\gamma y)=\gamma^*\mu_y(S_\infty)=\mu_y(S_\infty)=u(y)$ for $\gamma
\in\Gamma$, we have that $u$ is $\Gamma$-invariant. \par
Let $x\in\tilde{M}$ be fixed. It follows from, $ |\nabla B_\xi(y,x)|=|\cE|=1$ and 
Rauch'stheorem that we have
 $\exp(-DB_\xi(y,x))D(D-(n-1)b)
\leq\Delta\Theta^D\leq\exp(-D B_\xi(y,x))D(D-(n-1)a)$. This implies the result.
\end{proof} 
\begin{prop}\label{a40}
Suppose $\Gamma$ is nonelementary, i.e. has no abelian subgroup of finite index. Suppose that there are no nontrivial $\Gamma$-invariant positive-valued
superharmonic function on $\tilde{M}$. Then $(n-1)a\leq D\leq(n-1)b$. 
\end{prop}
\begin{proof}
Let $\Gamma x$ be the orbit of $x$ under $\Gamma$. Then the growth rate of 
the number of points of $\Gamma x$ in $\text{ball}(x,r)$ as $r$ increases
is bounded by $\text{vol}(\text{ball}(x,r))$. By the volume comparison theorem
we have $C_n\exp((n-1)br)\geq\text{vol}(\text{ball}(x,r))$, for some constant 
$C_n$ which depends only on dimension $n$. Therefore, when $s>(n-1)b$ we get 
$Z_\Gamma(x,s)<\infty$, which implies that $D\leq (n-1)b$.\par 
Next suppose that we have $D\leq (n-1)a$. Then by Proposition \ref{a39}, $u(x)$ is 
a $\Gamma$-invariant positive superharmonic and $\Delta u\leq D(D-(n-1)a) u$. 
It now follows from the hypothesis that $u$ is constant and that either $D=0$ or $D=(n-1)a$. 
However, since $\Gamma$ is nonelementary and $[\mu_y]$ is $\Gamma$-invariant, Proposition \ref{a34} 
 implies that $D\not= 0$. Hence, $D=(n-1)a$, and the result follows.
\end{proof}
\textnormal{The next proposition was originally proved by Sullivan \cite{Sul1} using a Borel-Cantelli 
type of argument. The proof is purely measure theoretic (see \cite{Nil}, \cite{Y2}). 
The proposition relates the ergodicity of $\G$ with the divergence of the
Poincar\'e series at} $D$. 
\begin{prop}[Sullivan]\label{a41}
Suppose that $\Gamma$ is nonelementary, discrete and torsion-free, and is divergent at $D$. Then $\G$ is ergodic with respect to $[\mu]^D$.  
\end{prop} 
\begin{prop}\label{a42}
Suppose $D=(n-1)a$ and there are no nontrivial positive
superharmonic functions on $M$. Then $\G$ is divergent.
\end{prop}
\begin{proof}
Fix a point $y\in\tilde{M}$. Let us assume the Poincar\'e series converges at $D$
(i.e. $\sum_{\g\in\G}\exp(-(n-1)a\dis(x,\g y))<\infty$). Then this series 
defines a nontrivial $\G$-invariant function on $\tilde{M}$. Let us 
denote this function by $h(x)$. Since 
$\exp(-D\dis(z,\g y))\le\exp(D\dis(x,z))\exp(-D\dis(x,\g y))$ for $z\in\tilde{M}$, it follows that for any
given  number $N>0$ there is a constant $C>0$ such that 
$\sum_{\g\in\G}\exp(-D\dis(z,\g y))\le C\sum_{\g\in\G}\exp(D\dis(x,\g y))$
for $\dis(z,x)\le N$. Hence the series converges uniformly on compact 
subsets of $\tilde{M}$. We will show that the convergence of the Poincar\'e
series at $(n-1)a$ implies existence of nontrivial positive superharmonic function on $M$.\par
Set $\dis_{\g y}(x):=\dis(x,\g y)$. First, we have 
\[ \Delta h(x)=\sum_{\g\in\G}\exp(-D\dis_{\g y}(x))D(D|\D\dis_{\g y}(x)|^2-\Delta\dis_{\g y}(x)).\] 
By Rauch'stheorem and $|\D\dis_{\g y}(x)|^2=1$ we get 
\[ \Delta h(x)\le\sum_{\g\in\G}\exp(-D\dis(x,\g y))D(D-(n-1)a),\] which 
implies $\Delta h\le 0$.
We consider the series $\sum_{\g\in\G}\log\tanh(\frac{(n-1)a\dis_{\g y}(x)}{2})$. It is 
easy to see that the convergence of this series on the set of points 
bounded away from $\G y$ follows from the convergence of the Poincar\'e 
series at $D=(n-1)a$. Denote this series by $-f$. Then by direct computation 
and Rauch's theorem we have $\Delta f(x)\le 0$ for 
$x\in\tilde{M}\backslash\G y$. Hence $1-\exp(-f(x))$ defines a nontrivial
positive $\G$-invariant superharmonic function on $\tilde{M}$. \par
Therefore, the convergence of the Poincar\'e series at $D=(n-1)a$ give raise 
to contradictions to our hypothesis, and the result follows. 
\end{proof}
\begin{cor}\label{a43}
Suppose $D=(n-1)a$ and there are no nontrivial positive
superharmonic functions on $M$. Then, $\G$ is ergodic with 
respect to $[\mu]^D$.
\end{cor}
\begin{proof}
The corollary follows from Proposition \ref{a41} and Proposition \ref{a42}.
\end{proof}
\begin{cor}\label{a44}
Let $M=\tilde{M}/\G$ be a topologically tame $3$-manifold with 
$-b^2\le\cK\le -1$ and $\Lambda(\G)=S_\infty$. If $D=2$, then $\G$ 
is divergent, hence ergodic with respect to $[\mu]^D$.
\end{cor}
\begin{proof}
The corollary follows from Theorem \ref{a16}, Proposition \ref{a42} and 
Corollary \ref{a43}.
\end{proof}
\begin{proof}[Proof of Theorem \ref{Main4}]
Under the hypothesis of Theorem \ref{Main4}, it follows from Proposition
\ref{a40} that $D\in [2,2b]$. That $\G$ is harmonically ergodic follows 
from Proposition \ref{a37}. If $D=2$, then
by Corollary \ref{a44} we have $\G$ is divergent.
\end{proof} 
\section{Part I of Theorems \ref{main-2} and \ref{main-1-2}}
Let $\G$ be a torsion-free discrete subgroup of $\iso(\widetilde M)$ with 
$D_\G=2$. We assume $\G$ is either convex-cocompact or 
$\Lambda_\G=S_\infty$, hausdorff-conservative and divergent.
\begin{prop}\label{eta-measure}
The measure $\mf{M}^2_{K_x}$ is finite and positive on all non-empty 
relative open subsets of 
$\Lambda_\G$, and $\mf{M}^2_{K_x}(A)=0$ if and only if $\mf{M}^2_{\eta_v}(A)=0$ for 
$A\subset\Lambda_\G\bs v(-\infty)$.
\end{prop}
\begin{proof}
First note that if we replace $\dis_{v,t}$ with
$\dis$ in the definition of $\eta_v$ we get a equivalent metric by 
Lemma 4 in \cite{Hem2}.\par
Let $x\in\widetilde M$ be any point. Denote $H_{v,x}$ the horosphere tangent
to $v(\infty)$ and passing through $x$. Take two vectors $U^\zeta, U^\xi$ 
in $S\widetilde M$ that are asymptotic to $v(-\infty)$ and passing through
$H_{v,x}$ with $U^\zeta(\infty)=\zeta$ and $U^\xi(\infty)=\xi$. Then 
there exists a positive constant $\alpha$ such that for
any unite tangent vectors $v^\zeta, v^\xi$ at $x$ which are asymptotic
to $\zeta$ and $\xi$ respectively, we have $\dis(g_t U^\zeta,g_t v^\zeta)
\le\alpha e^{-t}$ and $\dis(g_t U^\xi,g_t v^\xi)\le\alpha e^{-t}$, where
$g_t$ is the flow. This gives  
$\dis(g_\tau U^\zeta,g_\tau U^\xi)\le 2\alpha e^{-\tau}+1$ with 
$\tau=l_x(\zeta,\xi)$. On the other hand we also have 
$\beta^{-1}e^t\le\dis(g_t U^\zeta,g_t U^\xi)\le\beta e^{bt}$
for some positive constant $\beta$, which
gives $e^{-s}\ge\beta^{-1}$ and $e^{-s}\le\beta^{1/b}$ when  
$\dis(g_s U^\zeta, g_s U^\xi)=1$. Hence
$\frac{\beta^{-1}\beta^{1/b}}{2\alpha+1}e^{-s}\le e^{-\tau}\le\beta e^{-s}$.
Therefore $\eta_v$ and $K_x$ are equivalent on all points in $S_{v,x}$,
where $S_{v,x}$ is the shadow of $H_{v,x}$ cased from $v(-\infty)$. 
By compactness of $S_\infty$ there are $\{v_1,\dots,v_n\}\subset S\widetilde M$ 
such that
$\cup^n_1 S_{v_i,x}=S_\infty$. Since 
$0<\mf{M}^2_{\eta_{v_i}}(S_{v_i,x}\cap\Lambda_\G)<\infty$, we have 
$\mf{M}^2_{K_x}$ is positive and finite on $\Lambda_\G$. It follows
from Propositions \ref{a26}, \ref{prop-unique} and \ref{no-atom} the
measure $\mf{M}^2_{K_x}$ is positive on all relative open subsets.
Let $A\subset S_\infty\bs v(-\infty)$ be a $\mf{M}^2_{\eta_v}$-null set.
Let $\delta>0$.
Note that $\cup_{x\in\widetilde M}S_{v,x}=S_\infty\bs v(-\infty)$.
Hence there is $B\subset A$ with $B\subset S_{v,z}$ such that 
$\mf{M}^2_{K_z}(A\bs B)<\delta$. But 
$\mf{M}^2_{K_z}(B)\le c\mf{M}^2_{\eta_v}(B)$ for some $c>0$. By finiteness
we have $\mf{M}^2_{K_z}(A)<\delta$. Same argument holds for the rest of the 
proposition. 
\end{proof}
\begin{cor}\label{mu-eta}
The measures $\mu_x$ and $\mf{M}^2_{\eta_v}$ are absolutely continuous
with respect to each other. In-particular 
$\mf{M}^2_{\eta_v}$ is supported on $\Lambda_\G$.
\end{cor}
\begin{proof}
The result follows from Proposition \ref{a26}, \ref{prop-unique} and
Proposition \ref{eta-measure}.
\end{proof}
We use Mostow and Gehring's original idea to show the regularity of quasiconformal map \cite{Mos1}, \cite{Geh}.
This method was extended in \cite{Hem1}. We will follow their 
presentations, 
but with necessary generalizations that will allow us to prove
our theorems using results from previous sections.\par
Take the unite ball model of $\mathbb{H}^3$. Let $u$ be a unit tangent
vector at the origin. Let $\mathcal{O}_u$ be the unit circle on 
$\pr\mathbb{H}^3=S^2$ which is contained in the unique totally 
geodesic plane perpendicular
to $u$ and passing through the origin. Also denote the point $u(\infty)$ on $S^2$ by $\varsigma$. Then
for any pair $(p,\varsigma)\in\mathbb{B}_u:=\mathcal{O}_u\times\varsigma$ there 
is a unique 
semi-circle connecting them. The bundle of all these semi-circles is the 
upper hemisphere $\Omega_u$ of $S^2$. We denote this bundle space 
by $(\Omega_u,\pi_u,\mathbb{B}_u)$ 
where $\pi_u$ is the projection.\par
Let $\phi:S^2\longrightarrow S_\infty$ be a
quasi-conformal embedding conjugate $\G'$ to $\G$ 
under isomorphism $\chi:\G'\longrightarrow\G$, here
$\G'$ is a topologically tame, 
torsion-free, discrete subgroup of $\Psl(2,\mathbb{C})$ with 
$\Lambda_{\G'}=S^2$. And let $\psi$ be the inverse of $\phi$ 
when it is a quasi-conformal homeomorphism.\par
Let $\rho_u$ be the metric on $S^2\bs u(-\infty)$ which is defined same 
as $\eta_v$ with $v(-\infty)=\phi(u(-\infty))$. The hausdorff measure 
$\mf{M}^2_{\rho_u}$ on $S^2\bs u(-\infty)$ with
respect to $\rho_u$-metric is the usual Lebesgue measure. 
Hence there exists
a constant $\omega>0$ such that for all $\theta\in S^2\bs u(-\infty)$, 
we have $\mf{M}^2_{\rho_u}(B_{\rho_u}(\theta,r))=\omega r^2$. 
\begin{prop}[see \cite{Mos1}, \cite{Hem1}]\label{circles}
The measure $\phi^*\mf{M}^1_{\eta_v}$ is absolutely continuous with respect
to measure $\mf{M}^1_{\rho_u}$ on semi-circles. 
Here $\mf{M}^1_{\eta_v}$ and $\mf{M}^1_{\rho_u}$ are
$1$-dimensional hausdorff measures with respect to the $\eta_v$-metric and $\rho_u$-metric
respectively.
\end{prop}
\begin{proof}
Let $\mf{L}$ be the Lebesgue measure on $\mathbb{B}_u$. Then for 
all $P\in\mathbb{B}_u$ we have the following derivative
\[ \lambda(P):=\lim_{r\to 0}\frac{\mf{M}^2_{\eta_v}(\bar\phi\circ\pi^{-1}_u
(B_{\rho_u}(P,r)\cap\mathbb{B}_u))}{\mf{L}(B_{\rho_u}(P,r)\cap\mathbb{B}_u)}\]
exists and finite for $\mf{L}$-almost everywhere, see \cite{Fedder}.\par
Choose $P\in\mathbb{B}_u$ with $\lambda(P)<\infty$. 
For a semi-circle
$l:=\pi^{-1}_u(P)$, let $U_r(l)$ denote the $r$-neighborhood of $l$, then 
$\lim\sup_{r\to 0}\mf{M}^2_{\eta_v}(\bar\phi
(U_r(l)))/r<\infty$. For any compact $K\subset l$ with $\mf{M}^1_{\rho_u}(K)=0$, choose a number
$C>0$ with $\mf{M}^2_{\eta_v}(\bar\phi(U_r(l)))/r<C$. Let $\epsilon>0$
be given, by Besicovic's covering theorem there exists $\{\theta_1,\dots,\theta_k\}\subset K$ such that $kr<\epsilon$,
$K\subset\cup^k_1B_\rho(\theta_i,r)$ and any three of the balls $B_\rho(\theta_i,r)$
with distinct centers are disjoint.\par
Let $s_i:=\inf\{s>0|\bar\phi(B_\rho(\theta_i,r))\subset B_\eta(\bar\phi(\theta_i),s)\}$
and $\kappa>0$ (conformal constant) provided by Proposition \ref{qc-1}. Then we have 
$\bar\phi(K)\subset\cup^k_1 B_\eta(\bar\phi(\theta_i),s_i)$, 
$\bar\phi(S^2)\cap B_\eta(\bar\phi(\theta_i),s_i/\kappa)\subset\bar\phi(B_\rho(\theta_i,s_i))$. Since
$\G$ is hausdorff-conservative and by Proposition \ref{prop-1}, 
Corollary \ref{mu-eta},  
there exists $\alpha>0$ such that 
\begin{equation*}\begin{split}
\left(\sum^k_1 s_i\right)^2&\le k\sum^k_1 s^2_i\le k\kappa^2\alpha
\sum^k_1\mf{M}^2_{\eta_v}(\bar\phi(B_\rho(\theta_i,r)))\\
&\le 2\kappa^2\alpha k\mf{M}^2_{\eta_v}(\bar\phi(U_r(K)))\le 
2\kappa^2\alpha k\mf{M}^2_{\eta_v}(\bar\phi(U_r(l)))\\
&\le 2\kappa^2C\alpha(kr)\le\text{const}\quad\epsilon.
\end{split}\end{equation*}
Note the fact that any three of $\bar\phi(B_\rho(\theta_i,r))$ do not intersect
is used to bound $\sum^k_1\mf{M}^2_{\eta_v}(\bar\phi(B_\rho(\theta_i,r)))$ by
$2\mf{M}^2_{\eta_v}(\bar\phi(U_r(K)))$.\par
Therefore the result follows from the last inequality.
\end{proof}
The balls $B_\rho(\theta,r)$, $\theta\in S^2\bs u(-\infty)$, $r>0$ form a Vitali
relation for the Lebesgue measure $\mf{M}^2_{\rho_u}$. The following derivative
\[J(\theta):=\lim_{r\to 0}\frac{\mf{M}^2_{\eta_v}(\phi(B_\rho(\theta,r)))}
{\mf{M}^2_\rho(B_{\rho_u}(\theta,r))}\]
exists and finite for $\mf{M}^2_\rho$-almost 
every $\theta\in S^2\bs u(-\infty)$.
\begin{prop}\label{up-lower}
Let $\Lip_\phi$ be defined by 
$\Lip_\phi :\theta\longrightarrow\lim\sup_{r\to 0} r_\phi(\theta,r)/r$. 
Then $\Lip_\phi\in L^2_{\local} (S^2\bs u(-\infty),\mf{M}^2_{\rho_u})$. 
In-fact there exists a constant $k>0$ such that
\[\sqrt{J(\theta)}/k\le\lim\inf_{r\to 0}r_\phi(\theta,r)/r\le
\lim\sup_{r\to 0}r_\phi(\theta,r)/r\le k\sqrt{J(\theta)}.\]
\end{prop}
\begin{proof}
Let $\epsilon>0$. There is $r_\epsilon>0$ such that for any $r<r_\epsilon$
we have
\[\omega f(\theta)r^2/2\le\mf{M}^2_{\eta_v}(\bar\phi(B_\rho(\theta,r)))
\le (2\omega f(\theta)+\epsilon)r^2\]
where the fact that $\mf{M}^2_{\rho_u}$ is Lebesgue measure, 
i.e. $\mf{M}^2_{\rho_u}(B_{\rho_u}(\theta,r))=\omega r^2$ for some constant 
$\omega>0$ has been used. Since $\G$
is hausdorff-conservative and by Proposition \ref{prop-1},  
Corollary \ref{mu-eta}, 
there exists some constant $\alpha>0$ such that
\[(r_\phi(\theta,r)/\beta)^2/\alpha\le\mf{M}^2_\eta(\bar\phi(B_\rho(\theta,r)))
\le\alpha(r_\phi(\theta,r))^2.\]
Hence we have
\[\sqrt{(\omega/2\alpha f(\theta))}r\le r_\phi(\theta,r)\le 
\sqrt{\alpha(2f(\theta)\omega+\epsilon))}\beta r\]
and the result follows by letting $\epsilon\to 0$.
\end{proof}
\begin{lem}\label{loc-finite}
The image under $\phi$ of almost every 
semi-circle has locally finite $\mf{M}^1_{\eta_v}$-measure.
\end{lem}
\begin{proof}
Let $f:\Omega_u\longrightarrow\mathbb{B}_u\times [0,1]$ be a diffeomorphism
which maps $\pi^{-1}_u(P)$ over $P$ onto $P\times [0,1]$. For every
compact subset $C\subset\Omega_u$ we can find a positive number $\alpha$
such that
\begin{itemize}
\item For all $x\in f(C)$ the Jacobian of $f^{-1}$ at $x$ are $<\alpha$,
\item For all $P\in\mathbb{B}_u$ and $y\in\pi^{-1}_u(P)\cap C$ the local
dilations at $y$ of $f|_{\pi^{-1}_u(P)}$ are $<\alpha$.
\end{itemize}
Since $\bar\phi$ is a embedding, $\bar\phi(\Omega_u)$ is relative compact subset of $S_\infty\bs v(-\infty)$
and we have by Proposition \ref{up-lower},
$\int_{\Omega_u}\Lip^2_\phi\dif\mf{M}^2_{\rho_u}\le k^2\mf{M}^2_{\eta_v}(\bar\phi(\Omega_u))<\infty$, 
and H\"{o}lder inequality gives $\int_{\Omega_u}\Lip_\phi\dif\mf{M}^2_{\rho_u}<\infty$.
Hence 
\begin{equation*}\begin{split}
\int_{\mathbb{B}_u}\left(\int_{\pi^{-1}_u(P)\cap C}\Lip_\phi
\dif\mf{M}^1_{\eta_v}\right)\dif\mf{L}
&\le \alpha\int_{f(C)}\Lip_\phi\circ f^{-1}\dif\mf{L}\dif t\\
&\le\alpha^2\int_C\Lip_\phi\dif
\mf{M}^2_{\rho_u} <\infty
\end{split}\end{equation*}
where $\dif t$ is Lebesgue measure on $[0,1]$. 
Now by Proposition \ref{circles},  $\bar\phi$ is absolutely continuous on
$\pi^{-1}_u(P)$ therefore 
\[\mf{M}^1_{\eta_v}(\bar\phi(\pi^{-1}_u(P)\cap C))\le
\int_{\pi^{-1}_u(P)\cap C}\Lip_\phi\dif\mf{M}^1_{\rho_u}<\infty.\]
\end{proof}
Next we adapt the idea in \cite{Hem1} to prove the inequality part of 
Theorems \ref{main-2}.
\begin{proof}[Proof. Part I of Theorems \ref{main-2} and \ref{main-1-2}]\label{comp-D}
For Theorem \ref{main-1-2}, the inequality follows from Theorem \ref{Main4}
and Proposition \ref{lem-1}. 
Let $\G$ and $\G'$ be as in Theorem \ref{main-2} and satisfies those 
conditions.
Note that by Proposition \ref{lem-1}, $D_{\G'}=2$. Let $g$ be the Riemannian
metric of $\widetilde M$. Set $h=(D_\G/2)g$ as the new metric of $\widetilde M$.
The boundary space of $(\widetilde M,g)$ and $(\widetilde M,h)$ can be 
trivially identified, and 
$\eta_{(2/D_\G)v}=\eta^{D_\G/2}_v$. The critical exponent of $\G$ with
respect to $h$ is $2$, hence by Lemma \ref{loc-finite} there is a
non-trivial curve in $S_\infty\bs v(-\infty)$ with finite $D_\G/2$-dimensional
hausdorff measure with respect to $\eta_v$. However as noted before the 
curvature assumption $-b^2\le K\le -1$ of $g$ implies the $\eta_v$-metric 
is a distance on $S_\infty\bs v(-\infty)$, but the distance-hausdorff dimension
is $\ge 1$ for any non-trivial curves. Therefore we have $D_\G/2\ge 1$.
\end{proof}
\begin{lem}\label{nonsingular}
Let $\G'$ be a divergent, torsion-free discrete subgroup of 
$\Psl(2,\mathbb{C})$ with $\Lambda_{\G'}=S^2$ and $D_{\G'}=2$. 
Then the maps 
$\phi$ and $\psi$ are absolutely continuous with 
respect to $\sigma_y$ and $\mu_x$.
\end{lem}
\begin{proof}
By ergodicity of $\G$, $\G'$ and equivariance of $\bar\phi$, $\bar\psi$
and also Proposition \ref{mu-eta}, its suffices  
to show there exists a $A\subset S^2\bs u(-\infty)$ with 
$\mf{M}^2_{\rho_u}(A)>0$
such that the Radon-Nikodym derivative of $\bar\phi$ at every $x\in A$
with respect to $\mf{M}^2_{\rho_u}$ and $\mf{M}^2_{\eta_v}$ is non-zero. 
Using the
fact that $\eta_v$ is a distance function, it follows from 
Proposition \ref{circles}, for $\mf{L}$-almost all $P\in\mathbb{B}_u$ the length of
$\bar\phi(\pi^{-1}_u(P))>0$ is bounded by 
$\int_{\pi^{-1}_u(P)}\Lip_\phi\dif\mf{M}^1_{\rho_u}$. Hence if we set
$A:=\{x\in\Omega_u|\Lip_\phi(x)>0\}$, then for $\mf{L}$-almost all $P\in\mathbb{B}_u$,
$\mf{M}^1_{\rho_u}(\pi^{-1}_u(P)\cap A)>0$ which implies 
$\mf{M}^2_{\rho_u}(A)>0$.
Therefore the result follows from Proposition \ref{up-lower}.
\end{proof}
\section{Part II of Theorems \ref{main-2} and \ref{main-1-2}}
Let $\xi_1,\xi_2,\xi_3,\xi_4\in S_\infty$. The \emph{cross-ratio}
$|\xi_1,\xi_2,\xi_3,\xi_4|$ of these four points is defined as
\[|\xi_1,\xi_2,\xi_3,\xi_4|:=\frac{e^{-\beta_x(\xi_1,\xi_2)}e^{-\beta_x(\xi_3,\xi_4)}}
{e^{-\beta_x(\xi_1,\xi_3)}e^{-\beta_x(\xi_2,\xi_4)}}.\]
This definition is consistent with the hyperbolic space cross-ratio.\par
If $\G_1,\G_2$ are discrete subgroups of $\widetilde M$ such that both $\G_1,\G_2$
are divergent, and there exists a equivariant (under some group morphism $\chi$), 
nonsingular (with respect to $\mu_1,\mu_2$ Patterson-Sullivan measures
on $\Lambda_{\G_1}$ and $\Lambda_{\G_2}$ respectively), 
measurable map $f:\Lambda_{\G_1}\longrightarrow\Lambda_{\G_2}$. Then
\[\dif(f\times f)^*\Pi_2(\xi,\zeta)=e^{-D_{\G_2}\beta_y(f\xi,f\zeta)}
g(\xi)g(\zeta)\dif\mu_1(\xi)\dif\mu_2(\zeta)\]
where $g:=\frac{\dif f^*(\mu_2)}{\dif (\mu_1)}$, and $\Pi_i$ is
the measure defined in $\S 3$ through $\mu_i$. From
the properties of $f$, $(f\times f)^*\Pi_2$ is a constant $a>0$ multiple of 
$\Pi_1$. Hence 
$e^{D_{\G_2}\beta_y(f\xi,f\zeta)}g(\xi)g(\zeta)=ae^{D_{\G_1}\beta_x(\xi,\zeta)}$.
Therefore for $\mu_1$-almost everywhere we have
\[|f(\xi_1),f(\xi_2),f(\xi_3),f(\xi_4)|=|\xi_1,\xi_2,\xi_3,\xi_4|^{D_{\G_1}/D_{\G_2}}.\] 
This was the idea of Sullivan for the following lemma:
\begin{lem}\label{cr}
Let $\G_1,\G_2$ be discrete subgroups of $\iso(\widetilde M)$ with $D_{\G_1}=D_{\G_2}$ 
and $\G_1,\G_2$ are divergent. Suppose there exists a equivariant nonsingular
measurable map $f:\Lambda_{\G_1}\longrightarrow\Lambda_{\G_2}$ with respect
to Patterson-Sullivan measures space $(\Lambda_{\G_1},\mu_1)$ and 
$(\Lambda_{\G_2},\mu_2)$. Then $f$ preserves cross-ratio $\mu_1$-almost
everywhere. 
\end{lem}
For a finitely generated discrete subgroup $\G$ of $\Psl(2,\mathbb{C})$. 
The \emph{conservative set} of $\G$ on $S^2$ coincides 
with $\Lambda_{\G}$ 
up-to Lebesgue measure zero. The group $\G$ is called 
\emph{conservative} if and only if $\Lambda_{\G}$ has full 
Lebesgue measure. Since for a topologically tame $\G$, 
the hausdorff dimension of $\Lambda_\G$ is equal to $D_\G$, 
therefore we have the following:
\begin{prop}\label{hd-conservative}
Let $\G$ be a topologically tame, torsion-free discrete subgroup 
of $\Psl(2,\mathbb{C})$ with
conservative $\G$, then $\G$ is hausdorff-conservative.
\end{prop}
\begin{rem}\label{conser-hd} 
It is a conjecture that all finitely generated discrete subgroup $\G$ of
$\Psl(2,\mathbb{C})$ are topologically tame.
\end{rem}
Next we recall the statement of Sullivan's quasi-conformal stability 
for discrete subgroups of $\Psl(2,\mathbb{C})$.
\begin{thm}[Sullivan \cite{Sul}]\label{sul-conf}
Let $\G$ be a discrete subgroup of $\Psl(2,\mathbb{C})$. Then $\G$ is
quasi-conformally stable (i.e. if $f$ is a quasi-conformal automorphism of $S^2$ 
with $f\G f^{-1}\subset\Psl(2,\mathbb{C})$, then $f$ is a M\"{o}bius transformation) 
if and only if $\G$ is conservative.
\end{thm} 
\begin{cor}\label{sul-rig}
Let $N=\mathbb{H}^3/\G$ be a complete hyperbolic $3$-manifold for
a conservative $\G$. Then $N$ is quasi-isometrically stable, 
i.e. If there is a quasi-isometric homeomorphism $h:N\longrightarrow M$
to a hyperbolic manifold $M$, then $N$ is isometric to $M$.
\end{cor}
\begin{proof}[Proof. Theorem \ref{main-1-2} part II]
By Theorem \ref{Main4}, $\G$ is divergent for $D_\G=2$. 
From Proposition \ref{lem-1}, $\G'$ is also divergent and $D_{\G'}=2$. 
Lemma \ref{nonsingular} then implies $f$ is absolutely 
continuous with respect to $\sigma_y$ and $\mu_x$. 
Hence by Lemma \ref{cr}, $f$
preserves cross ratio $\sigma_y$-everywhere. By Proposition \ref{prop-1}, $\Lambda_{\G'}=S^2$
and since $\sigma_y$ is non-zero constant multiple 
of Lebesgue measure, we can modify $f$ on the Lebesgue measure
null subset of $S^2$ to a map which is cross ration preserving on $S^2$.
We denote the new map also by $f$. 
By Bourdon's theorem \cite{Bourdon},
$f$ extends into the space as a isometry, i.e. $\mathbb{H}^3$ and
$\widetilde M$ are isometric. Hence the result follows from 
Theorem \ref{sul-conf}.
\end{proof}
\begin{proof}[Proof. Theorems \ref{main-2} part II]
Here $f$ embeds $S^2$ into $S_\infty$. If we suppose 
$D_\G=D_{\G'}=2$, then by using same argument as the proof of 
Theorem \ref{main-1-2}, $f$
extends to a isometric embedding of $\mathbb{H}^3$ into $\widetilde M$ by \cite{Bourdon}.
Since $f(S^2)$ is a $\Lambda_\G$-invariant closed subset of 
$S_\infty$, by Proposition \ref{prop-1}, $f(S^2)=\Lambda_\G$. Hence the boundary space of
the isometric embedded image of $\mathbb{H}^3$ coincides with $\Lambda_\G$, 
therefore the result follows.
\end{proof} 
\begin{proof}[Proof. Corollary \ref{cor-1}]
This follows from Propositions \ref{qc-1}, \ref{lem-1}, and Theorem \ref{main-1-2}.
\end{proof}
\bibliography{Ref}

\end{document}